\documentclass[12pt]{amsart}
\usepackage{pinlabel}
\usepackage{amsfonts,amsmath,amscd,xcolor,graphicx,amssymb}
\usepackage{caption}
\usepackage{epstopdf}
\usepackage{a4wide}

\theoremstyle{plain}

\newtheorem*{Lemma 2.5}{Lemma 2.5}

\theoremstyle{definition}

\theoremstyle{remark}

\makeatletter
\@namedef{subjclassname@2020}{\textup{2020} Mathematics Subject Classification}
\makeatother

\begin{document}

\title[On the Invalidity of Lemma 2.5]{On the Invalidity of Lemma 2.5 in our previous work on the Powell Conjecture}

\author[S. Cho]{Sangbum Cho}
\thanks{The first-named author is supported by the National Research Foundation of Korea(NRF) grant funded by the Korea government(MSIT) (NRF-2021R1F1A1060603 and RS-2024-00456645).}
\address{Department of Mathematics Education, Hanyang University, Seoul 04763, Korea, and School of Mathematics, Korea Institute for Advanced Study, Seoul 02455, Korea}
\email{scho@hanyang.ac.kr}

\author[Y. Koda]{Yuya Koda}
\thanks{The second-named author is supported by JSPS KAKENHI Grant Numbers JP23K20791, JP23H05437, and 24K06744}
\address{Department of Mathematics, Hiyoshi Campus, Keio University, Yokohama 223-8521, Japan, and
International Institute for Sustainability with Knotted Chiral Meta Matter (WPI-SKCM$^2$), Hiroshima University, Higashi-Hiroshima 739-8526, Japan}
\email{koda@keio.jp}

\author[J. H. Lee]{Jung Hoon Lee}
\thanks{The third-named author is supported by the National Research Foundation of Korea(NRF) grant funded by the Korea government(MSIT) (RS-2023-00275419).}
\address{Department of Mathematics and Institute of Pure and Applied Mathematics, Jeonbuk National University, Jeonju 54896, Korea}
\email{junghoon@jbnu.ac.kr}

\author[N. Sekino]{Nozomu Sekino}
\email{nnsekino@gmail.com}

\date{\today}

\begin{abstract}
In our previous version entitled ``The reducing sphere complexes for the 3-sphere are connected: a proof of the Powell Conjecture", we claimed to prove the Powell Conjecture, which states that the Goeritz group of the genus-$g$ Heegaard splitting of the 3-sphere is finitely generated for any non-negative integer $g$. However, we have found a critical error in the proof of Lemma 2.5 in that version. 
In this note, we prove that the statement of Lemma 2.5 does not hold in general. 
This invalidates a key step in our argument and leaves the proof of the Powell Conjecture incomplete. 
Consequently, the Powell Conjecture remains an open problem in the case of $g \geq 4$.
\end{abstract}

\maketitle

The Powell Conjecture offers specific finite generating sets for the Goeritz group of the Heegaard splitting of the $3$-sphere. 
For the genus-$2$ Heegaard splitting, this was established by Goeritz~\cite{Goe33} and Scharlemann~\cite{Sch04}. 
The case of genus $3$ is also valid by Freedman-Scharlemann \cite{Freedman-Scharlemann}, as well as Cho-Koda-Lee~\cite{CKL24}. 
In our previous version, titled ``The reducing sphere complexes for the 3-sphere are connected: a proof of the Powell Conjecture", 
we claimed  to establish the connectivity of the reducing sphere complexes of the Heegaard splittings of the 3-sphere, and thereby prove the Powell Conjecture for any genus due to Zupan~\cite{Zupan}. 
However, we have found a critical flaw in the proof of Lemma 2.5 in that version, which we recall below, prompting us to investigate potential counterexamples.
As a result, we have constructed explicit examples showing that the statement of Lemma 2.5 is false in general. 
The following sections present these counterexamples for genus $g \geq 3$. 

\begin{Lemma 2.5}
Let $(V, W; \Sigma)$ be a genus-$g$ Heegaard splitting of $S^3$, where $g \geq 2$.
Let $D_1, \ldots, D_n$ be pairwise disjoint, pairwise non-parallel essential disks in $V$.
Let $E$ be a non-separating disk in $W$ disjoint from $\bigcup_{i=1}^n D_i$.
Then there exists a reducing sphere $P$ for $\Sigma$ disjoint from $(\bigcup_{i=1}^n D_i) \cup E$ such that $\bigcup_{i=1}^n D_i$ and $E$ are in different components of $S^3 \setminus P$.
\end{Lemma 2.5}

\section{The case of $g = 3$}

\begin{figure}[h]
\centering
\vspace{1em}

\includegraphics[width=14cm]{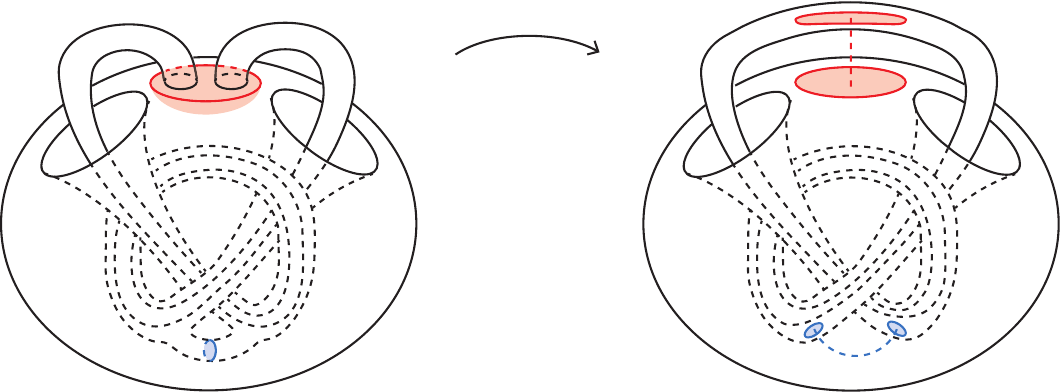}
\begin{picture}(400,0)(0,0)
\put(174,156){Compress}
\put(168,125){along $D \cup E$}
\put(76,110){\footnotesize \color{red} $D$}
\put(77,16.5){\footnotesize \color{blue} $E$}
\put(125,42){$V$}
\put(145,20){$W$}
\put(40,-5){$\Sigma := \partial V = \partial W$}

\put(365,42){$V_0$}
\put(298,-5){$\Sigma_0 := \partial V_0$}
\end{picture}
\vspace{0.5em}
\captionof{figure}{The genus-$3$ surface $\Sigma$.}
\label{fig:counter_example_genus-3}
\end{figure}

The surface $\Sigma$ in Figure~\ref{fig:counter_example_genus-3} is a genus-$3$ Heegaard surface that cuts $S^3$ into two handlebodies $V$ and $W$. 
Consider the separating essential disk $D$ in $V$ and  the non-separating essential disk $E$ in $W$ shown in the figure.
If we compress $\Sigma$ along $D$ and $E$, we obtain two (parallel) knotted tori.
One of the two tori, denoted by $\Sigma_0$, contains two scars coming from $E$ and one scar coming from $D$.
As in Figure~\ref{fig:counter_example_genus-3}, $\Sigma_0$ bounds the exterior of a knotted solid torus, which we denote by $V_0$, and of course $V_0$ contains no essential disks.

We claim that there is no reducing spheres for $\Sigma$ that separates $D$ and $E$.
Suppose, for contradiction, that we have such a reducing sphere $P$. 
Then, if we compress $\Sigma$ along the disk $P \cap V$, the resulting two surfaces are both Heegaard surfaces of $S^3$, and so $P \cap V$ cuts $V$ into two ``unknotted'' handlebodies in $S^3$. 
Since the torus component of the surface obtained by compressing $\Sigma$ along $D$ is the boundary of a regular neighborhood of a trefoil, 
the separating circle $C = P \cap \Sigma$ is not isotopic to $ \partial D$ in $\Sigma$. 
Therefore, the circle $C = P \cap \Sigma$ must lie in $\Sigma_0$, and so the disk $P \cap V$ lies in $V_0$. 
We denote by $V_1$ one of the two handlebodies containing the disk $D$.

Now, since $V_0$ contains no essential disks, the disk $P \cap V$ is inessential in $V_0$, implying that $C$ is also inessential in $\Sigma_0$.
Let $F$ be the disk in $\Sigma_0$ bounded by $C$. 
We have two possibilities, either $F$ contains only a single scar coming from $D$ or $F$ contains only the two scars coming from $E$. 
The former is impossible because if so, $\partial D$ and $\partial F = C$ is isotopic in $\Sigma$, a contradiction. 
The latter case implies that $V_1$ is not a handlebody, since in this case $V_1$ is isotopic to the union of $V$ and a regular neighborhood of $E$, which is again a contradiction. 

Thus, we conclude that there is no reducing spheres for $\Sigma$ separating $D$ and $E$.

\vspace{1em}

We remark that this is an example of the weak reducing pair $D - E$ of a separating disk $D$ and a non-separating disk $E$ which does not allow a reducing sphere separating $D$ and $E$. 
But, if both of $D$ and $E$ in the weak reducing pair are all non-separating, then we can always find a reducing sphere separating 
$D$ and $E$.

\section{The case of $g \geq 3$}

\begin{figure}[h]
\centering
\includegraphics[width=14cm]{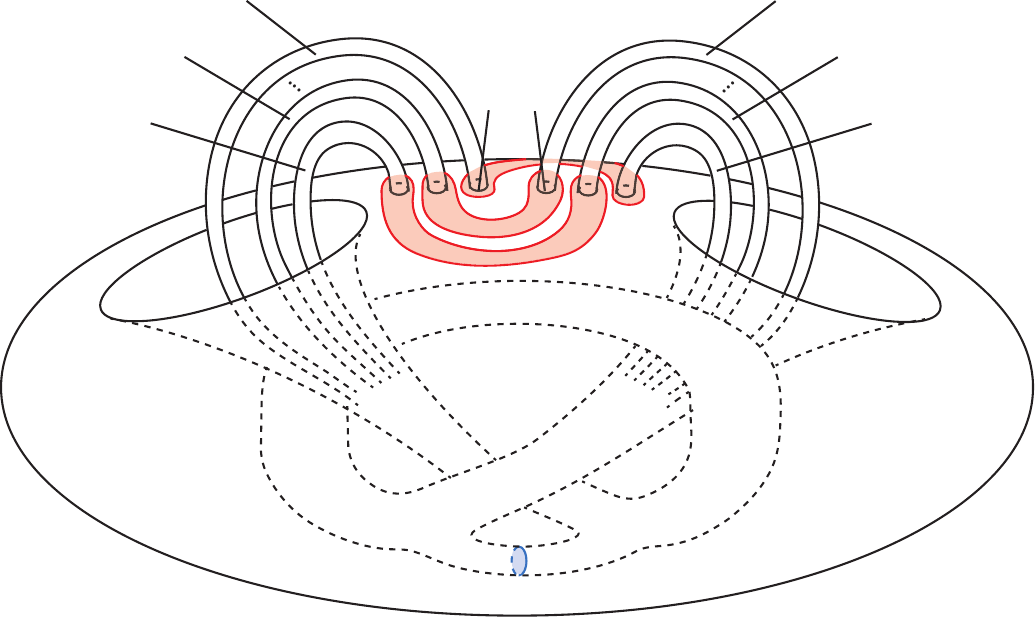}
\begin{picture}(400,0)(0,0)
\put(47,200){$h_1$}
\put(57,230){$h_2$}
\put(72,250){$h_{g-2}$}
\put(342,200){$h_1$}
\put(328,230){$h_2$}
\put(305,250){$h_{g-2}$}
\put(153,169){\scriptsize $l_1$}
\put(169,169){\scriptsize $l_2$}
\put(182,213){\scriptsize $l_{g-2}$}
\put(201,213){\scriptsize $r_{g-2}$}
\put(225,170){\scriptsize $r_2$}
\put(241,166){\scriptsize $r_1$}
\put(160,155){\color{red} \scriptsize $D_1$}
\put(175,163){\color{red} \scriptsize $D_2$}
\put(189,193){\color{red} \scriptsize $D_{g-2}$}
\put(197,18){\color{blue} $E$}
\put(340,60){\Large $V$}
\put(380,30){\Large $W$}
\put(160,-5){$\Sigma := \partial V = \partial W$}

\end{picture}
\vspace{0.5em}
\captionof{figure}{The genus-$g$ surface $\Sigma$.}
\label{fig:counter_example_general}
\end{figure}

We will generalize the example of genus $3$ almost directly.
The example for $g=3$ that we have seen in the previous section is  exactly the special case of the following construction. 
In the case of $g=3$, we had a single handle attached on the ``two-holed sphere'' whose feet lie in the disk bounded by $\partial D$.
Now, for $g \geq 3$, we attach $g-2$ copies $h_1, h_2, \ldots, h_{g-2}$ of the handle to have a genus-$g$ Heegaard surface $\Sigma$ as in Figure~\ref{fig:counter_example_general}.
For each $i \in \{1, 2, \ldots, g-2\}$, we denoted by $l_i$ and $r_i$ the two feet of the handle $h_i$ as in the figure.
Let $V$ and $W$ be the genus-$g$ handlebodies cut off by $\Sigma$, and again we may assume that the non-separating essential disk $E$ is contained in $W$.

\vspace{1em}

Next, we fix any bijection $\phi:\{1, 2, \ldots, g-2\} \rightarrow \{1, 2, \ldots, g-2\}$, 
and choose $g-2$ essential disks $D_1, D_2, \ldots, D_{g-2}$ in $V$ satisfying that
\begin{itemize}
\item the boundary circles $\partial D_j$ lies in the ``two-holed sphere'', and
\item for each $j \in \{1, 2, \ldots, g-2\}$, the circle $\partial D_j$ bounds a disk in the two-holed sphere containing only $l_j$ and $r_{\phi(j)}$ among all the feet of the handles $h_1, h_2, \ldots, h_{g-2}$.
\end{itemize}
In the figure, we choose a bijection $\phi$ such that $\phi(1) = 2$, $\phi(2) = g-2$, $\ldots, \phi(g-2) = 1$, for example.

Note that, if $g=3$, we have only a single disk $D = D_1$, which is necessarily separating, but if $g \geq 4$, some of the disks $D_j$ might be separating and others non-separating, depending on the choice of the bijection $\phi$.
Actually, we can choose a suitable $\phi$ easily so that the resulting disks $D_j$ are all separating or all non-separating.
If we compress $\Sigma$ along $D_1 \cup D_2 \cup \cdots \cup D_{g-2}$ and $E$, we obtain a bunch of knotted tori.
One of them, denoted by $\Sigma_0$, contains two scars coming from $E$ and exactly one scar coming from each $D_j$. 
Each of the other components $T_1, T_2, \ldots, T_k$ is the boundary of a regular neighborhood of a trefoil or of a non-trivial satellite of a trefoil, where the number $k$ ($1 \leq k \leq g-2$) depends on the choice of $\phi$. 
(We observe that $k = g-2$ if and only if $\phi$ is the identity map.)
We also note that each of $T_1, T_2, \ldots, T_k$ contains no scars coming from $E$. 
As in the case of $g = 3$, $\Sigma_0$ bounds the complement of a knotted solid torus, which we denote by $V_0$, and again $V_0$ contains no essential disks.

We claim that there is no reducing spheres for $\Sigma$ that separates $D_1 \cup D_2 \cup \cdots \cup D_{g-2}$ and $E$.
Suppose, for contradiction, we have such a reducing sphere $P$. 
Then, by compressing $\Sigma$ along the disk $P \cap V$, we obtain two surfaces, which are both Heegaard surfaces of $S^3$, and so $P \cap V$ cuts $V$ into two ``unknotted'' handlebodies in $S^3$. 
If the circle $C = P \cap \Sigma$ lies in one of $T_1, T_2, \ldots, T_k$, say $T_1$, then, by compressing $\Sigma$ along $P \cap V$, we obtain a ``knotted" torus component isotopic to $T_1$, which is a contradiction. 
Therefore, the circle $C$ lies in $\Sigma_0$, and so the disk $P \cap V$ lies in $V_0$. 
We denote by $V_1$ one of the two handlebodies containing $D_1 \cup D_2 \cup \cdots \cup D_{g-2}$.

Now, since $V_0$ contains no essential disks, the disk $P \cap V$ is inessential in $V_0$, implying that $C$ is also inessential in $\Sigma_0$.
Let $F$ be the disk in $\Sigma_0$ bounded by $C$.
We note that the two disks $F$ and $P \cap V$ are parallel in $V$, and they have the common boundary circle $C$. 
We have two possibilities, either $F$ contains only the two scars coming from $E$ or $F$ contains only the scars coming from $D_1, D_2, \ldots, D_{g-2}$.
The former case implies that $V_1$ is not a handlebody, since in this case $V_1$ is isotopic to the union of $V$ and a regular neighborhood of $E$, that is a contradiction. 
The latter case implies that $V_1$ is a ``knotted'' handlebody (of genus $g-2$), since in this case $V_1$ is the manifold which is the union of $g-2$ ``solid'' handles (of $h_1, h_2, \ldots, h_{g-2}$) and the $3$-ball in $V$ bounded by $F \cup (P \cap V)$, a contradiction again.

Thus, we conclude that there is no reducing spheres for $\Sigma$ separating $D_1 \cup D_2 \cup \cdots \cup D_{g-2}$ and $E$.

\bibliographystyle{amsplain}


\end{document}